\documentclass[12pt,a4paper]{article}
\usepackage{epsfig}
\usepackage{amsmath}

\usepackage{latexsym}
\usepackage{amssymb}
\usepackage{paralist}
\usepackage{amsfonts}

\newtheorem{thm}{Theorem}[section]

\newtheorem{prop}[thm]{Proposition}
\newtheorem{lem}[thm]{Lemma}
\newtheorem{cor}[thm]{Corollary}

\newcommand\Q[1]{$#1$ \quad}

\font\msbm=msbm10 at 12pt
\newcommand{\ZZ}{\mbox{\msbm Z}}

\newcommand{\FF}{\mbox{\msbm F}}

\newcommand{\FG}{F \left[ G^{\ast} \right]}
\newcommand{\KG}{K \left[ G^{\ast} \right]}
\newcommand{\G}{G^{\ast}}

\begin{document} 

\title{On the Existence of\\ Hermitian Self-Dual\\ Extended Abelian Group Codes}

\author{Lilibeth Dicuangco\footnote{Mathematics Department, University of the Philippines, Diliman, Quezon 
City, 1101 Philippines, {\tt ldicuangco@math.upd.edu.ph}}, \: \setcounter{footnote}{6} Pieter 
Moree\footnote{Max-Planck-Institut, Vivatsgasse 7, D-53111 Bonn, Germany, {\tt moree@mpim-bonn.mpg.de}},  
\: \setcounter{footnote}{3} Patrick Sol\'e\footnote{CNRS, I3S ESSI, BP 145, Route des Colles, 06 903 
Sophia Antipolis, France, {\tt sole@essi.fr}}}

\date{} 
\maketitle 
\date{} 
\begin{abstract}
Split group codes are a class of group algebra codes over an abelian group. They were introduced by Ding, 
Kohel and Ling in \cite{cunsheng} as a generalization of the cyclic duadic codes. For a prime power $q$ 
and an abelian group $G$ of order $n$ such that gcd$(n,q)=1$, consider the group algebra 
$\FF_{q^{2}}[G^{*}]$ of $\FF_{q^{2}}$ over the dual group $G^{*}$ of $G$. We prove that every ideal code 
in $\FF_{q^{2}}[G^{*}]$ whose extended code is Hermitian self-dual is a split group code. We characterize 
the orders of finite abelian groups $G$ for which an ideal code of $\FF_{q^{2}}[G^{*}]$ whose extension is 
Hermitian self-dual exists and 
derive asymptotic estimates for the number of non-isomorphic abelian groups 
with this property.
\end{abstract}



\section{Introduction}

Binary duadic codes were first introduced in 1984 by Leon, Masley and Pless \cite{leon et al} as a 
generalization of quadratic residue codes. Smid \cite{smidthesis} generalized these results further by 
defining duadic codes over arbitrary finite fields in terms of a \emph{splitting} of the length of the 
code.  The Q-codes of Pless \cite{qcodes} are then duadic codes over $\FF_{4}$ in this setting.

Quadratic residue codes have also been generalized in a different direction (i.e., see \cite{ward}). In 
this approach, quadratic residue codes are defined as ideals of abelian group algebras, which is a  
generalization of cyclic codes. Rushanan \cite{rushanan} proceeded to define duadic codes in this setting.

In \cite{cunsheng}, Ding, Kohel and Ling defined split group codes as ideals of abelian group algebras. 
Their construction makes use of a \emph{splitting} of the abelian group. Under this definition, split 
group codes are seen as a generalization of duadic codes.

In this paper, we consider the finite field $F = \FF_{q^2}$ and an abelian group $G$ of order $n$ such 
that gcd$(n,q)=1$. Following the treatment in \cite{cunsheng}, we work with the dual group $G^{\ast}$ of 
$G$ and consider the group algebra $F\left[ G^{\ast} \right]$. We prove 
that every ideal code in $\FG$ whose extension by a suitable parity-check is Hermitian self-dual is a 
split group code (Corollary \ref{Cor exte2}).
We then give sufficient and necessary conditions on the order of the group $G$ for the existence of 
Hermitian self-dual extended  ideal codes (Theorem \ref{Thm concon2}). We conclude the
paper by deriving asymptotic estimates on $HSD(x)$, the number of non-isomorphic abelian groups of order 
less than $x$ for which a Hermitian self-dual extended ideal code exists (Theorem \ref{twee}).

\section{Preliminaries}

All of the results in this section are taken from \cite{cunsheng}. In general, these results work for any 
finite field but since we will be dealing with Hermitian duality, we restrict our study to finite fields 
of square order.

Let $R$ be a finite commutative ring with unity. Let $G$ be its underlying finite abelian group written 
additively. Denote the order and exponent of $G$ by $n$ and $m$, respectively. Let $q$ be a power of a 
prime $p_{1}$ such that gcd$(n,q) = 1$ or equivalently, gcd$(m,q) = 1$. Let $F = \FF_{q^{2}}$. Let $K$ be 
the smallest extension of $F$ containing all the $m$-th roots of unity.
 
Let $G^{\ast}$ be the set of all characters of $G$ into $K$. The groups $G$ and $\G$ are isomorphic. Let 
$K \left[ G^{\ast} \right]$ be the group algebra of $K$ over $G^{\ast}$. The elements of $\KG$ are the 
sums $\sum_{\psi \in \G} a_{\psi}\psi$ where the $a_{\psi}$'s are elements of $K$. An ideal $I$ of $\KG$ 
is called an \emph{ideal code}.

The dimension of the commutative group algebra $\KG$ over $K$ is $n$. This group algebra $\KG$ contains a 
subgroup isomorphic to $G^{\ast}$. For any character $\psi$ in $\G$, we also denote by $\psi$ the 
corresponding element in $\KG$. 

If $x \in G$ and $f = \sum_{\psi \in \G} a_{\psi}\psi \in \KG$, define $f(x) = \sum_{\psi \in \G} 
a_{\psi}\psi(x)$.  Thus we can view the elements of $\KG$ as functions from $G$ to $K$.

\subsection{Ideal codes and Idempotent Generators}

An element $e$ of a ring is called an idempotent if $e^{2} = e$. An idempotent is called primitive if for 
every other idempotent $f$, either $ef = e$ or $ef = 0$.

\begin{prop} {\rm (\cite{cunsheng}.)} The primitive idempotents of $\KG$ are the elements
\begin{equation*}
e_{x} = \frac{1}{n} \sum_{\psi \in \G} \psi(x)^{-1}\psi,
\end{equation*}
for each element $x$ of $G$.
\end{prop}

\begin{prop} {\rm (\cite{cunsheng}.)} The ring $\KG$ decomposes as a direct sum $\oplus_{x \in G} K 
e_{x}$. If $f \in \KG$, then $f$ has the form
\begin{equation*}
f = \sum_{x \in G} f(x) e_{x}.
\end{equation*}
Every idempotent $e$ in $\KG$ can be uniquely written in the form
\begin{equation*}
e = \sum_{x \in X} e_{x},
\end{equation*}
for some non-empty subset $X$ of $G$.
\end{prop}

\noindent Let $X$ be a non-empty subset of $G$ and define the ideal
\begin{equation*}
I_{X} = \{ f \in \KG \mid f(x) = 0 \; \mbox{ for all } x \in X \}.
\end{equation*}

\begin{cor} {\rm (\cite{cunsheng}.)}
For every ideal $I$ in $\KG$, there is a unique proper subset $X$ of $G$ such that $I = I_{X}$ and $I$ is 
generated by the idempotent $e = \sum_{x \notin X} e_{x}$.
\end{cor}

\subsection{Split Group Codes}

Let $s$ be an element of $R$. Consider the endomorphism of $G$ given by $\tau_{s} : x \longrightarrow sx$. 
This induces a map $\mu_{s}$ on $\G$  given by $\mu_{s}(\psi) = \psi\circ\tau_{s}$ for each element $\psi$ 
of $\G$. This extends to a map on $\KG$, also denoted by $\mu_{s}$, defined by $\mu_{s}(f) = 
f\circ\tau_{s}$ for all $f \in \KG$. That is, for $f \in \KG$, $\mu_{s}(f)(x) = f(sx)$ for every $x \in 
G$.

A \emph{splitting} of $G$ over $Z$ is a triple $(Z,X_{0},X_{1})$ which gives a partition $G = Z \cup X_{0} 
\cup X_{1}$ such that there exists an invertible element $s$ of $R$ with $\tau_{s}(X_{0}) = X_{1}$ and 
$\tau_{s}(X_{1}) = X_{0}$. Under these conditions, $s$ is said to \emph{split} the triple 
$(Z,X_{0},X_{1})$. In addition, an invertible element $r$ of $R$ is said to \emph{stabilize} the splitting 
if $\tau_{r}(X_{0}) = X_{0}$ and $\tau_{r}(X_{1}) = X_{1}$.

Given a splitting $(Z,X_{0},X_{1})$, let $C_{0}(K)$ be the ideal $I_{X_{0}}$ over $K$ and let $C_{1}(K)$ 
be the ideal $I_{X_{1}}$ over $K$. The ideal $C_{0}(K)$ is defined as the \emph{split group code} 
associated to the splitting, and the ideal $C_{1}(K)$ is called the \emph{conjugate split group code}. The 
following notations are used to denote some special subcodes: $C_{0}^{Z}(K) = I_{Z \cup X_{0}}$,  
$C_{1}^{Z}(K) = I_{Z \cup X_{1}}$ and $C_{Z}(K) = I_{X_{0} \cup X_{1}}$.

Let $s$ be an invertible element of $R$. The element $s$ is said to \emph{split} the group code $C_{0}(K)$ 
if $\mu_{s}(C_{0}(K)) = C_{1}(K)$ and $\mu_{s}(C_{1}(K)) = C_{0}(K)$, while $s$ is said to 
\emph{stabilize} the code $C_{0}(K)$ if $\mu_{s}(C_{0}(K)) = C_{0}(K)$ and $\mu_{s}(C_{1}(K)) = C_{1}(K)$.

\begin{prop} {\rm (\cite{cunsheng}.)}
Let $s$ be a unit in $R$. A split group code $C_{0}(K)$ is split or stabilized by $s$ if and only if $s$ 
splits or stabilizes $(Z,X_{0},X_{1})$, respectively.
\label{s}
\end{prop}

\begin{thm} {\rm (\cite{cunsheng}.)}
Let $(Z,X_{0},X_{1})$ be a splitting. Let $C_{0}(K)$ be the split group code associated to this splitting. 
Then the following hold:
\begin{enumerate}
\item The codes $C_{0}(K)$ and $C_{1}(K)$ are generated by the idempotents
\begin{equation*}
e = \sum_{x \notin X_{0}} e_{x} \; \mbox{ and } \; f = \sum_{x \notin X_{1}} e_{x}.
\end{equation*}
The codes $C_{0}^{Z}(K)$, $C_{1}^{Z}(K)$ and $C_{Z}(K)$ are generated by 
\begin{equation*}
\sum_{x \in X_{1}} e_{x}, \; \sum_{x \in X_{0}} e_{x}, \; \mbox{ and } \; \sum_{z \in Z} e_{z}.
\end{equation*}
\item If the splitting is given by $s$, then $\mu_{s}$ induces an equivalence of $C_{0}(K)$ with its 
conjugate $C_{1}(K)$, and of the subcode $C_{0}^{Z}(K)$ with $C_{1}^{Z}(K)$.
\item $\KG$ decomposes as a direct sum $C_{Z}(K) \oplus C_{0}^{Z}(K) \oplus C_{1}^{Z}(K)$.
\end{enumerate}
\label{equi}
\end{thm}

\begin{cor} {\rm (\cite{cunsheng}.)}
The codes $C_{0}(K)$ and $C_{1}(K)$ have dimension $(n + |Z|)/2$. The subcodes $C_{0}^{Z}(K)$ and 
$C_{1}^{Z}(K)$ have dimension $(n - |Z|)/2$. The subcode $C_{Z}(K)$ has dimension $|Z|$.
\end{cor}

\subsection{Split Group Codes Over $\FF_{q^{2}}$}

In the previous section, the split group codes are defined over the field $K$ which is assumed to contain 
all the $m$-th roots of unity. In this paper, we want our split group codes to be defined over the 
subfield $F =\FF_{q^{2}}$ without requiring $\FF_{q^{2}}$ to contain any $m$-th roots of unity.
In this section, we present sufficient and necessary conditions given in \cite{cunsheng} for split group 
codes to be defined over a subfield of $K$.

Let $V = V(K)$ be a vector subspace of $K^{n} = \KG$. Define $V(F) = V(K) \cap F^{n}$. Clearly, 
$\mbox{dim}_{F}(V(F)) \leq \mbox{dim}_{K}(V(K))$. If equality holds then we say that $V$ is \emph{defined 
over the field $F$}.

If the vector subspace $C_{0}(K)$ of $K^{n} = \KG$ is defined over $F$, we simply write $C_{0}$ for the 
subcode $C_{0}(F) = C_{0}(K) \cap F^{n}$ in $F^{n} = \FG$. In this case, we call $C_{0}$ \emph{the split 
group code over $F$}. Similarly, we write $C_{1}$, $C_{0}^{Z}$, $C_{1}^{Z}$ and $C_{Z}$ for the other 
codes defined over $F$.

Note that $(m,q)=1$ by assumption, so the integer $-q$ as an element of the finite ring $R$ is invertible 
and $\tau_{q^2}$ is a well-defined automorphism of $G$. The action of the group generated by 
$\tau_{q^{2}}$ on the elements of $G$ partitions $G$ into disjoint orbits. These $\langle \tau_{q^{2}} 
\rangle$-orbits play the same role as the cyclotomic cosets for the cyclic codes.

\begin{prop} {\rm (\cite{cunsheng}.)}
The idempotents of $\FG$ are those $e$ in $\KG$ of the form
\begin{equation*}
e = \sum_{x \in Y} e_{x},
\end{equation*}
where $Y$ is a union of $\langle \tau_{q^{2}} \rangle$-orbits in $G$. An idempotent $e$ in $\FG$ is 
primitive if and only if $Y = \langle \tau_{q^{2}} \rangle x$ for some $x \in G$.
\end{prop}

\begin{cor} {\rm (\cite{cunsheng}.)}
Let $\{ e_{1}, e_{2}, \ldots, e_{r} \}$ be the set of all primitive idempotents of $\FG$. Then every 
nonzero ideal $I$ of $\FG$ is generated by $e = \sum_{i \in T} e_{i}$  where $T$ is a non-empty subset of 
$\{1, 2,\ldots, r \}$.
\end{cor}

\begin{thm} {\rm (\cite{cunsheng}.)}
Let $I$ be an ideal in $\KG$. Then the following conditions are equivalent.
\begin{enumerate}
\item The ideal $I$ is defined over $F$.
\item The set $X = \{ x \in G \: | \: f(x) = 0 \; \mbox{ for all } f \in I  \}$ is a union of $\langle 
\tau_{q^{2}} \rangle$-orbits.
\item The idempotent generator of $I$ lies in $\FG$. 
\end{enumerate}
\label{descent}
\end{thm}

\begin{cor} {\rm (\cite{cunsheng}.)}
If $C_{0}(K)$ is defined over $F$ then so is $C_{1}(K)$. Moreover $\FG$ has the decomposition
\begin{equation*}
\FG = C_{Z}(F) \oplus C_{0}^{Z}(F) \oplus C_{1}^{Z}(F). 
\end{equation*}
If $s$ gives the splitting, then $\mu_{s}$ gives an equivalence of $C_{0}$ and $C_{1}$ and of  $C_{0}^{Z}$ 
and $C_{1}^{Z}$.
\end{cor}

\begin{thm} {\rm (\cite{cunsheng}.)}
Let $C$ be a code in $\KG$.  The block length, dimension, and minimum distance are well-defined invariants 
of $C$, independent of the field over which $C$ is defined.
\end{thm}

\section{Hermitian Duality and Extended Ideal Codes in $\FG$}

In this section, we present results concerning the Hermitian orthogonality of ideal codes in $\FG$. We 
prove that every Hermitian self-orthogonal ideal code in $\FG$ is a subcode of a split group code for some 
splitting of $G$ given by $-q$. We then proceed to define an extension of  
ideal codes
and determine conditions for the extended code to be Hermitian self-dual. We also give necessary and 
sufficient conditions on the order of the abelian group $G$ for the existence of Hermitian self-dual 
extended ideal codes.

\subsection{Hermitian Orthogonality of Ideal Codes in $\FG$}

Let $f = \sum_{\psi \in \G} a_{\psi} \psi$ and $g = \sum_{\psi \in \G} b_{\psi} \psi$ be elements of 
$\FG$. The \emph{Hermitian inner product} between $f$ and $g$ is defined as $\langle f, g \rangle_{H} = 
\sum_{\psi \in \G} a_{\psi}b_{\psi}^{q}$.

Let $C$ be a code in $\FG$. The Hermitian dual of $C$ is the ideal $C^{\bot_{H}} = \{ f \in \FG \; | \; 
\langle f, g \rangle_{H} = 0 \; \mbox{ for all } g \in C \}$. The code $C$ is said to be \emph{Hermitian 
self-orthogonal} if $C \subseteq C^{\bot_{H}}$ and is \emph{Hermitian self-dual} if $C = C^{\bot_{H}}$.

Theorem \ref{Thm HermSplit} states the main result of this section. It is a generalization of Proposition 
4.4 in \cite{DMS} to split group codes. We first prove some basic results concerning the Hermitian duals 
of split group codes.

\begin{prop}
Let $f = \sum_{\psi \in \G} a_{\psi} \psi$ and $g = \sum_{\psi \in \G} b_{\psi} \psi$ be elements of 
$\FG$. Then the Hermitian inner product of $f$ and $g$ is
\begin{equation*}
\langle f, \: g \rangle_{H} = \frac{1}{n}\sum_{x \in G}f(x)g(-q^{-1}x)^{q}
\end{equation*}
\label{hermdot}
\end{prop}
{\it Proof}. Note that 
$\mu_{-q^{-1}}(g)(x) = g(-q^{-1}x)  =  \sum_{\psi \in \G}b_{\psi}\psi(-q^{-1}x)
                         =  \sum_{\psi \in \G}b_{\psi} \psi(q^{-1}x)^{-1}.$
Thus
 $(\mu_{-q^{-1}}(g)(x))^{q}  =  \sum_{\psi \in \G}b_{\psi}^{q} \psi(q^{-1}x)^{-q}
 						    =  \sum_{\psi \in \G}b_{\psi}^{q} \psi(x)^{-1},$			
or equivalently, $(\mu_{-q^{-1}}(g))^{q} = \sum_{\psi \in \G}b_{\psi}^{q} \psi^{-1}$.

Define $f \ast g = f(\mu_{-q^{-1}}(g))^{q}$.
Then the coefficient of the trivial character of $f \ast g$ is $\sum_{\psi \in \G} a_{\psi}b_{\psi}^{q}$, 
which is $\langle f, g \rangle_{H}$.
Expanding $f \ast g$ in terms of its idempotent decomposition, we get
\begin{eqnarray*}
 f \ast g & = & \sum_{x \in G}(f(\mu_{-q^{-1}}(g))^{q})(x)e_{x} \\
          & = & \sum_{x \in G}f(x)(\mu_{-q^{-1}}(g)(x))^{q}e_{x} \\
          & = & \sum_{x \in G}f(x)(g(-q^{-1}x))^{q}e_{x} \\
          & = & \frac{1}{n} \sum_{\psi \in \G} \sum_{x \in G} f(x)(g(-q^{-1}x))^{q} \psi^{-1}(x)\psi.
\end{eqnarray*}
Using this expansion, the coefficient of the trivial character of $f \ast g$ is $\frac{1}{n} \sum_{x \in 
G}f(x)(g(-q^{-1}x))^{q}$. The result follows.
\hfill \Q{\Box}

\begin{prop}
Let $C$ be an ideal in $\KG$ which is defined over $F$. Suppose $C = I_{X}$ for some non-empty subset $X$ 
of $G$. Then $C^{\bot_{H}} = I_{X'}$ where $X' = G \setminus \tau_{-q}(X)$.
\label{hermdual}
\end{prop}
{\it Proof}. {}From Theorem \ref{descent}, $X$ is a union of $\langle \tau_{q^{2}} \rangle$-orbits.
Let $X' = G \setminus \tau_{-q}(X)$. Let $f \in I_{X'}$ and $g \in C = I_{X}$. By Proposition 
\ref{hermdot}, $\langle f, g \rangle_{H} = \frac{1}{n}\sum_{x \in G}f(x)(g(-q^{-1}x))^{q}$.  Since $X$ and 
$\tau_{-q}(X)$  are unions of $\langle \tau_{q^{2}} \rangle$-orbits and clearly $x$ and $q^{2}x$ belong to 
the same $\tau_{q^{2}}$-orbit, it follows that
$-q^{-1}x \in X$ if and only if $(-q)^{2}(-q^{-1})x \in X$ 
                 if and only if $-qx \in X$
                 if and only if $(-q)(-q)x \in \tau_{-q}(X)$
				 if and only if  $x \in \tau_{-q}(X)$. 
Thus $g(-q^{-1}x) = 0$ for all $x \in \tau_{-q}(X)$. Since $f \in I_{X'}$, $f(x) = 0$ for every $x \in X' 
= G \setminus \tau_{-q}(X)$. Therefore $\sum_{x \in G} f(x)(g(-q^{-1}x))^{q} = 0$, implying that $\langle 
f, g \rangle_{H} = 0$. Thus $I_{X'} \subseteq C^{\bot_{H}}$. Comparing dimensions, we get $C^{\bot_{H}} = 
I_{X'}$.
\hfill \Q{\Box}

\vskip .1in
\noindent {\tt Remark}. We note that for any subset $X$ of $G$ which is a union of $\langle \tau_{q^{2}} 
\rangle$-orbits of $G$, we have $\tau_{-q^{-1}}(X) = \tau_{(-q)^{2}}(\tau_{-q^{-1}}(X)) = \tau_{-q}(X)$.

\begin{prop}
Let $(Z,X_{0},X_{1})$ be a splitting of $G$ over $Z$ where $Z$, $X_{0}$ and $X_{1}$ are unions of $\langle 
\tau_{q^{2}} \rangle$-orbits. The ring element $-q$ splits or stabilizes $C_{0}$ if and only if $-q^{-1}$ 
splits or stabilizes $C_{0}$, respectively.
\label{s2}
\end{prop}
{\it Proof}. Using Proposition \ref{s}, we need only show that $-q$ splits or stabilizes $(Z,X_{0},X_{1})$ 
if and only if  $-q^{-1}$  splits or stabilizes $(Z,X_{0},X_{1})$, respectively. {}From the remark above, 
$\tau_{-q^{-1}}(X_{0}) = \tau_{-q}(X_{0})$ and $\tau_{-q^{-1}}(X_{1}) = \tau_{-q}(X_{1})$. The result 
follows.
\hfill \Q{\Box}

\begin{prop}
Let $(Z,X_{0},X_{1})$ be a splitting of $G$ over $Z$. Assume that $Z$, $X_{0}$ and $X_{1}$ are unions of 
$\langle \tau_{q^{2}} \rangle$-orbits.  Suppose $Z$ is stabilized by $\tau_{-q}$. Then $C_{Z}^{\bot_{H}} = 
C_{0}^{Z} \oplus C_{1}^{Z}$. If $-q$ splits $C_{0}$, then $C_{0}^{\bot_{H}} = C_{0}^{Z}$. If $-q$ 
stabilizes $C_{0}$, then $C_{0}^{\bot_{H}} = C_{1}^{Z}$.
\label{OddSplit}
\end{prop}
{\it Proof}. Note that $\tau_{-q^{-1}}(Z) = \tau_{(-q)^{2}}(\tau_{-q^{-1}}(Z)) = \tau_{-q}(Z)$ and by 
assumption, $\tau_{-q}(Z) = Z$.
Let $f \in C_{Z} = I _{X_{0} \cup X_{1}}$ and let $g \in C_{0}^{Z} \oplus C_{1}^{Z}$. Then
\begin{equation*}
\langle f, g \rangle_{H} = \frac{1}{n} \sum_{x \in G}f(x)(g(-q^{-1}x))^{q} = 0,
\end{equation*}
since $f(x) = 0$ for all $x \in X_{0} \cup X_{1}$ and $g(-q^{-1}x) = 0$ for all $x \in Z$. Thus $C_{0}^{Z} 
\oplus C_{1}^{Z} \subseteq C_{Z}^{\bot_{H}}$. Comparing dimensions, we have $C_{Z}^{\bot_{H}} = C_{0}^{Z} 
\oplus C_{1}^{Z}$.

Suppose $-q$ splits $C_{0}$. Let $f \in C_{0}$ and $g \in C_{1}^{Z}$. Note that $\tau_{-q^{-1}}(Z) = Z$. 
By assumption, $\tau_{-q}(X_{1}) = X_{0}$, or equivalently, $\tau_{-q^{-1}}(X_{1}) = X_{0}$. Clearly $f(x) 
= 0$ for all $x \in X_{0}$ and $g(-q^{-1}x) = 0$ for all $x \in Z \cup X_{1}$. Thus 
\begin{equation*}
\langle f, g \rangle_{H} = \frac{1}{n} \sum_{x \in G}f(x)(g(-q^{-1}x))^{q} = 0,
\end{equation*}
implying that $C_{0}^{Z} \subseteq C_{0}^{\bot_{H}}$. Comparing dimensions, we have $C_{0}^{\bot_{H}} = 
C_{0}^{Z}$.

Suppose $-q$ stabilizes $C_{0}$. Then $\tau_{-q^{-1}}(X_{1}) = \tau_{-q}(X_{1}) = X_{1}$. Let $f \in 
C_{0}$ and $g \in C_{1}^{Z}$. Clearly $f(x) = 0$ for all $x \in X_{0}$ and  $g(-q^{-1}x) = 0$ for all $x 
\in Z \cup X_{1}$. Thus 
\begin{equation*}
\langle f, g \rangle_{H} = \frac{1}{n} \sum_{x \in G}f(x)(g(-q^{-1}x))^{q} = 0,
\end{equation*}
implying that $C_{1}^{Z} \subseteq C_{0}^{\bot_{H}}$. Comparing dimensions, we have $C_{1}^{\bot_{H}} = 
C_{0}^{Z}$.
\hfill \Q{\Box}

\vskip .2in

\noindent We now prove the main result of this section.

\begin{thm}
Let $C = I_{X}$ be an ideal in $\FG$. Then $C$ is Hermitian self-orthogonal if and only if $C = C_{0}^{Z}$ 
for some splitting $(Z, X_{0}, X_{1})$ of $G$ which is split by $-q$ (that is, $C$ is a subcode of a split 
group code which is split by $-q$).
\label{Thm HermSplit}
\end{thm}
{\it Proof}.

\noindent $(\Longleftarrow)$ Suppose $(Z, X_{0}, X_{1})$ is a splitting of $G$ which is split by $-q$. Let 
$C = C_{0}^{Z}$. By Proposition \ref{OddSplit}, $C^{\bot_{H}} = (C_{0}^{Z})^{\bot_{H}} = C_{0} \supseteq 
C_{0}^{Z} = C$.

\noindent $(\Longrightarrow)$ Let $C^{\bot_{H}} = I_{X'}$. By Proposition \ref{hermdual}, $X' = G 
\setminus \tau_{-q}(X)$. By assumption, $C \subseteq C^{\bot_{H}}$. Thus $X' \subseteq X$, or $G \setminus 
\tau_{-q}(X) \subseteq X$. Note that both $X$ and $X'$ are unions of $\langle \tau_{q^{2}} 
\rangle$-orbits. Write $X = Z \cup X_{0}$ with
$\tau_{-q}(Z) = Z \; \mbox{ and } \; Z \cap X_{0} = \emptyset$.
Since $\tau_{-q}$ either fixes a $\langle \tau_{q^{2}} \rangle$-orbit or it sends it to another $\langle 
\tau_{q^{2}} \rangle$-orbit, we are actually choosing $Z$ as the union of $\langle \tau_{q^{2}} 
\rangle$-orbits contained in $X$ which are fixed by $\tau_{-q}$ and $X_{0}$ as the complement of $Z$ in 
$X$, so that clearly $\tau_{-q}(X_{0}) \cap X_{0} = \emptyset$.

We first show that neither $Z$ nor $X_{0}$ is empty. If $0 \notin X$ then $0 \notin \tau_{-q}(X)$, which 
implies that $0 \in G \setminus \tau_{-q}X$. But $G \setminus \tau_{-q}(X) \subseteq X$, implying that $0 
\in X$, a contradiction. Thus $0 \in X$, and so by our choice of partition of $X$, $0 \in Z$ proving that 
$Z \not= \emptyset$. If $X_{0} = \emptyset$ then $X = Z$ and $\tau_{-q}(X) = X$. Thus $X' = G \setminus 
\tau_{-q}(X) = G \setminus X$ which implies that $X' \cap X = \emptyset$, a contradiction since $X' 
\subseteq X$. Thus $X_{0} \not= \emptyset$.

Note that $X' = G \setminus \tau_{-q}(X) = G \setminus Z \cup \tau_{-q}(X_{0})$. Since $X_{0} \cap Z = 
X_{0} \cap \tau_{-q}(X_{0}) = \emptyset$, we have $X_{0} \subseteq X'$. Consider $\tau_{-q}(X') = G 
\setminus X$. Since $X' \subseteq X$, it follows that $\tau_{-q}(X') \cap X' = \emptyset$
and $\tau_{-q}$ does not fix any $\langle \tau_{q^{2}} \rangle$-orbit contained in $X'$, implying that $X' 
\subseteq X_{0}$. Hence $X' = X_{0}$.
Let $X_{1} = \tau_{-q}(X')$. Then $(Z, X_{0}, X_{1})$ gives a splitting of $G$ such that $\tau_{-q}(Z) = 
Z$, $\tau_{-q}(X_{0}) = X_{1}$ and $\tau_{-q}(X_{1}) = X_{0}$  and $C = I_{X} = I_{Z \cup X_{0}} = 
C_{0}^{Z}$.
\hfill \Q{\Box}

\subsection{Extensions of Ideal Codes in $\FG$}

Using the terminology of split group codes, duadic codes are easily seen to be split group codes for 
splittings of $G$ over $Z = \{0 \}$ where $G$ is cyclic (see Example III.1 of \cite{cunsheng}). In 
\cite{DMS}, we defined an extension for an odd-like duadic code and gave a sufficient condition for the 
extended code to be Hermitian self-dual. In this section, we consider split group codes for the abelian 
group $G$ with splittings over $Z = \{0 \}$ and derive analogous results regarding Hermitian self-duality 
of the extended split group codes.

Let the order $n$ of the abelian group $G$ be odd. Consider the equation
\begin{equation}
\frac{1}{n} + \gamma^{q+1} = 0,
\label{gam}
\end{equation}
which is solvable in $\FF_{q^{2}}$.
Let $\gamma$ be a solution to (\ref{gam}). For each $f \in \FG$, define $\widetilde{f} = (f, -\gamma f(0)) 
\in \FG \times F$. If $C$ is a code in $\FG$ then the extended code $\widetilde{C}$ is defined as the 
subspace
\begin{equation*}
\widetilde{C} = \{ \widetilde{f} = (f, -\gamma f(0)) \mid f \in C \} \subseteq \FG \times F.
\end{equation*}

 \begin{prop}
 Let $(Z = \{0 \}, X_{0}, X_{1})$ be a splitting of $G$ where $Z$, $X_{0}$ and $X_{1}$ are unions of 
$\langle \tau_{q^{2}} \rangle$-orbits. Let  $C_{0}$ be the corresponding split group code defined over $F 
= \FF_{q^{2}}$.
\begin{enumerate}
\item The extended codes $\widetilde{C_{0}}$ and $\widetilde{C_{1}}$ are equivalent. 
\item If $-q$ splits $C_{0}$, then $\widetilde{C_{0}}^{\bot_{H}} = \widetilde{C_{0}}$ and 
$\widetilde{C_{1}}^{\bot_{H}} = \widetilde{C_{1}}$.
\item If $-q$ stabilizes $C_{0}$, then $\widetilde{C_{0}}^{\bot_{H}} = \widetilde{C_{1}}$ and 
$\widetilde{C_{1}}^{\bot_{H}} = \widetilde{C_{0}}$.
\end{enumerate}
\label{Thm exte2}
 \end{prop}
{\it Proof}. The equivalence of $\widetilde{C_{0}}$ and $\widetilde{C_{1}}$ is an immediate consequence of 
Theorem \ref{equi}.

Suppose $-q$ splits $C_{0}$. Let $f$, $g$ be elements of $C_{0}$. {}From Proposition \ref{s2}, $-q^{-1}$ 
also splits $C_{0}$. Thus $\mu_{-q^{-1}}(g) \in C_{1}$. Note that $-q^{-1}x \in X_{0} \iff x \in 
\tau_{-q^{-1}}(X_{0}) = X_{1}$. It follows that $f(x) = 0$ for all $x \in X_{0}$ and $g(-q^{-1}x) = 0$ for 
all $x \in X_{1}$. Hence using Proposition \ref{hermdot}, we get
\begin{eqnarray*}
\langle f, g \rangle_{H} & = & \frac{1}{n} f(0) g(0)^{q} = -\gamma^{q+1} f(0) g(0)^{q}.
\end{eqnarray*}
Thus $\langle \widetilde{f}, \widetilde{g} \rangle_{H} = 0$ and $\widetilde{C_{0}}^{\bot_{H}} = 
\widetilde{C_{0}}$. By a similar argument, it can be shown that $\widetilde{C_{1}}^{\bot_{H}} = 
\widetilde{C_{1}}$.

Suppose $-q$ stabilizes $C_{0}$. Let $f \in C_{0}$ and let $g \in C_{1}$. The element $-q^{-1}$ also 
stabilizes $C_{0}$ and $\mu_{-q^{-1}}(g) \in C_{0}$. Again, $-q^{-1}x \in X_{1} \iff x \in 
\tau_{-q^{-1}}(X_{1}) = X_{1}$. So $f(x) = 0 $ for all $x \in X_{0}$ and $g(-q^{-1}x) = 0$ for all $x \in 
X_{1}$. Thus
\begin{eqnarray*}
\langle f, g \rangle_{H} & = & \frac{1}{n} f(0) g(0)^{q} = -\gamma^{q+1} f(0) g(0)^{q},
\end{eqnarray*}
and so $\langle \widetilde{f}, \widetilde{g} \rangle_{H} = 0$ and $\widetilde{C_{0}}^{\bot_{H}} = 
\widetilde{C_{1}}$. Similarly, $\widetilde{C_{1}}^{\bot_{H}} = \widetilde{C_{0}}$.
\hfill \Q{\Box}

\begin{cor}
Let $C = I_{X}$ be a group code defined over $F$. The extended code $\widetilde{C}^{\bot_{H}}$ is 
Hermitian self-dual if and only if $C$ is a split group code for some splitting $(Z = \{0 \}, X_{0}, 
X_{1})$ of $G$ by $-q$.
\label{Cor exte2}
\end{cor}
{\it Proof}. 

\noindent $(\Longleftarrow)$ This follows directly from the previous theorem.

\noindent $(\Longrightarrow)$ Since $\widetilde{C}^{\bot_{H}}$ is Hermitian self-dual, the dimension of 
$C$ is $\frac{n+1}{2}$ and so $C$ cannot be Hermitian self-orthogonal. This fact combined with the 
assumption that $\widetilde{C}^{\bot_{H}}$ is Hermitian self-dual implies that $0 \notin X$.
Let $C_{e} = I_{X \cup \{ 0 \} }$. This subcode $C_{e}$ is Hermitian self-orthogonal and has dimension 
$\frac{n-1}{2}$. By Theorem \ref{Thm HermSplit}, $C_{e}$ is a subcode of a split group code which is split 
by $-q$, that is, $C_{e} = C_{0}^{Z}$ for some splitting $(Z, X_{0}, X_{1})$ of $G$ by $-q$. Since dim 
$C_{e} = \frac{n-1}{2}$ and dim $C_{0}^{Z} = \frac{n- |Z|}{2}$, it follows that $Z = \{0 \}$. Hence $X = 
X_{0}$ and $C$ is a split group code of $G$ which is split by $-q$.
\hfill \Q{\Box}

\subsection{Existence of Hermitian Self-dual Extended Ideal Codes}

In view of Theorem \ref{Thm HermSplit} and Corollary \ref{Cor exte2}, it is natural to ask under what 
conditions we obtain splittings over $Z = \{ 0 \}$ of an abelian group $G$ by $-q$. Such conditions would 
guarantee existence of Hermitian self-orthogonal codes and Hermitian self-dual extended codes in $\FG$. We 
remark that the results in this section are generalizations of results on extended duadic codes in 
\cite{DMS}.

Define $ord_{r}(q)$ to be the smallest positive integer $t$ such that $q^{t} \equiv 1 \;(\mbox{mod } r)$. 
If $l$ is a positive odd integer relatively prime to $q$ then $l$ is said to be split by $-q$ over 
$\FF_{q^{2}}$ if and only if the set $X = \{1, 2, \ldots, l \}$ has a partition $X = X_{0} \cup X_{1}$ 
such that $(-q)X_{0} = X_{1}$ and $(-q)X_{1} = X_{0}$, where the multiplication is read modulo $l$.

\begin{prop} {\rm (\cite{DMS}.)}
Let $l$ be a positive odd integer which is relatively prime to $q$. The integer $l$ has a splitting by 
$-q$ if and only if $ord_{r}(q) \not\equiv 2 \;(\mbox{mod } 4)$ for every prime $r$ dividing $l$.
\label{duadic}
\end{prop}

\begin{thm}
Let $G$ be an abelian group of order $n$. The group $G$ has a splitting over $Z = \{ 0 \}$ given by $-q$ 
if and only if $ord_{r}(q) \not\equiv 2 \;(\mbox{mod } 4)$ for every prime $r$ dividing $n$.
\label{spli}
\end{thm}
{\it Proof}.
The abelian group $G$ is isomorphic to a unique product of cyclic groups of the form
\begin{equation*}
\ZZ_{m_{1}} \times \ZZ_{m_{2}} \times \ldots \times \ZZ_{m_{s}},
\end{equation*}
where $m_{i}$ divides $m_{i+1}$ for $i = 1, 2, \ldots, s-1$, and $m_{s} = m$ where $m$ denotes the 
exponent of $G$.

If each summand $\ZZ_{m_{i}}$ has a splitting over $Z = \{ 0 \}$ given by $-q$ then $G$ also has a 
splitting over $Z = \{ 0 \}$ given by $-q$. Indeed if $(\{ 0 \}, X_{0}^{(i)}, X_{1}^{(i)})$ is a splitting 
by $-q$ of $\ZZ_{m_{i}}$ for each $i = 1, 2, \ldots, s$ then letting
\begin{eqnarray*}
X_{t} & = & \hspace{.1in} X_{t}^{(1)} \times \ZZ_{m_{2}} \times \ZZ_{m_{3}} \times \cdots \times 
\ZZ_{m_{s}}\\
	  &   & \cup \; \{ 0 \} \times X_{t}^{(2)} \times \ZZ_{m_{3}} \times \cdots \times \ZZ_{m_{s}}\\
	  &   & \cup \; \{ 0 \} \times \{ 0 \} \times X_{t}^{(3)} \times \ZZ_{m_{4}} \times \cdots \times 
\ZZ_{m_{s}} \\
	  &   & \hspace{1in} \vdots \\
	  &   & \cup \; \{ 0 \} \times \{ 0 \} \times \{ 0 \} \times \cdots \times \{ 0 \} \times X_{t}^{(s)}
\end{eqnarray*}
for $t = 0, 1$, $(\{ 0 \}, X_{0}, X_{1})$ gives a splitting of $G$ by $-q$.
Conversely, suppose that $G$ has a splitting over $Z = \{ 0 \}$ given by $-q$. Let $\ZZ_{i} = \{ 0 \} 
\times \{ 0 \} \times \ldots \times \ZZ_{m_{i}} \times \{ 0 \} \times \ldots \{ 0 \}$ be the subgroup of 
$G$ isomorphic to $\ZZ_{m_{i}}$. If $(Z = \{0 \}, X_{0}, X_{1})$ gives a splitting for $G$ by $-q$, then 
$(Z = \{0 \}, \ZZ_{i} \cap X_{0}, \ZZ_{i} \cap X_{1})$ gives a splitting for $\ZZ_{i}$ given by $-q$.
Hence using Proposition \ref{duadic},
$G$ has a splitting over $Z = \{ 0 \}$ given by $-q$
if and only if each summand $\ZZ_{m_{i}}$ has a splitting over $Z = \{ 0 \}$ given by $-q$
if and only if $m_{i}$ is split by $-q$ for all $i = 1, 2, \ldots , m$
if and only if  $ord_{r}(q) \not\equiv 2 \;(\mbox{mod } 4)$ for every prime $r$ dividing $m_{i}$ for all 
$i = 1, 2, \ldots , s$
if and only if  $ord_{r}(q) \not\equiv 2 \;(\mbox{mod } 4)$ for every prime $r$ dividing $m$.
But the primes dividing $m$ are precisely the primes dividing $n$.
Thus $G$ has a splitting over $Z = \{ 0 \}$ given by $-q$ if and only if $ord_{r}(q) \not\equiv 2 
\;(\mbox{mod } 4)$ for every prime $r$ dividing $n$.
\hfill \Q{\Box}

\vskip .1in

\noindent \emph{Example}:  Let $G = \ZZ_{3} \times \ZZ_{9}$ and $\FF_{4^{2}} = \FF_{16}$. Note that 
$ord_{3}(4) = 1$ and by Theorem \ref{spli} the abelian group $G$ has a partition which is split by $-4$. 
The cyclic groups $\ZZ_{3}$ and $\ZZ_{9}$ have splittings by the multiplier $\mu_{-4}$ given by $(\{0\}, 
A_{1}, A_{2})$ and $(\{0\}, B_{1} \cup B_{3}, B_{2} \cup B_{6})$, respectively, where $A_{i}$ is the 
$16$-cyclotomic coset modulo $3$ containing $i$ and $B_{j}$ is the $16$-cyclotomic coset modulo $9$ which 
contains $j$. Define $C_{(i,j)}$ as the orbit of $\tau_{16}$ in $G$ containing $(i,j)$. Letting $X_{0} = 
C_{(1,0)} \cup C_{(1,1)}  \cup C_{(1,2)} \cup C_{(1,3)}  \cup C_{(1,6)}  \cup C_{(0,1)}  \cup C_{(0,3)}$ 
and $X_{1} = C_{(2,0)} \cup C_{(2,1)}  \cup C_{(2,2)} \cup C_{(2,3)}  \cup C_{(2,6)}  \cup C_{(0,2)}  \cup 
C_{(0,6)}$, the set $(\{ (0,0) \}, X_{0}, X_{1} )$ gives a splitting of $G$ by $-4$. Notice that this 
partition can be obtained from the splittings of $\ZZ_{3}$ and $\ZZ_{9}$ as described in the proof.
\hfill \Q{\Box}

\vskip .1in

\noindent We remark that $ord_{r}(q) \not\equiv 2 \;(\mbox{mod } 4)$ means that either $ord_{r}(q)$ is odd 
or $ord_{r}(q)$ is doubly even. It can easily be verified that $ord_{r}(q)$ is doubly even if and only if 
$ord_{r}(q^{2})$ is even. Thus Theorem \ref{spli} can be restated as:

\begin{thm}
 Let $G$ be an abelian group of order $n$. The group $G$ has a splitting over $Z = \{ 0 \}$ given by $-q$ 
if and only if  for every prime $r$ dividing $n$,  either $ord_{r}(q)$ is odd or $ord_{r}(q^{2})$ is even.
\label{Thm spli2}
\end{thm}

\noindent Thus we get the following condition for the existence of extended ideal codes of $\FG$ which are 
Hermitian self-dual.

\begin{thm}
Let $G$ be an abelian group of order $n$.
An ideal code of $\FG$ whose extension is Hermitian self-dual exists if and only if for every prime $r$ 
dividing $n$,  either $ord_{r}(q)$ is odd or $ord_{r}(q^{2})$ is even.
\label{Thm concon2}
\end{thm}
{\it Proof}. This is a direct consequence of Corollary \ref{Cor exte2} and Theorem \ref{Thm spli2}.
\hfill \Q{\Box}

\vskip .1in

\noindent We note that the same result as the preceding theorem was obtained by Mart\'{i}nez-P\'{e}rez and 
Willems in \cite{conchita} for ideal codes in a group algebra over any finite group.

\section{Counting Hermitian self-dual extended abelian group codes}
Theorem \ref{Thm concon2} raises the question of counting the number of non-isomorphic
abelian groups of order $\le x$ for which an ideal code of $F[G^*]$ whose
extension is Hermitian self-dual exists. An estimate for this quantity
is provided by Theorem \ref{twee}.

Let $q=p_1^t$ be a prime power and let ${\cal P}_q$ be the set of primes 
$r\ne p_{1}$ for which ord$_{r}(q)$ is odd or ord$_{r}(q^2)$ is even. 
Let ${\cal P}_q(x)$
be the associated counting function. The primes $r$ not counted, that
is the primes $r$ such that ord$_r(q)\equiv 2({\rm mod~}4)$  or the
prime $r=p_1$ can 
be shown, see \cite{DMS}, to have a natural density $\delta(q)$ that
is given by the following formula (with $\lambda$ the exponent of 2 in the
factorisation of $t$):
\begin{equation*}
\delta(q)=\delta(p_1^{t})=\left\{ 
\begin{array}{ll}
7/24                & \mbox{ if }  p_1=2 \mbox{ and }  \lambda=0 ; \\
1/3                 & \mbox{ if } p_1=2 \mbox{ and } \lambda=1 ; \\
2^{-\lambda -1}/3   & \mbox{ if } p_1=2 \mbox{ and } \lambda\ge 2 ;\\
2^{-\lambda}/3      & \mbox{ if } p_1\ne 2.
\end{array}
\right.
\end{equation*}
It can be proved, see \cite[Lemma A.3]{DMS}, that 
\begin{equation}
\label{estimate}
{\cal P}_q(x)=(1-\delta(q)){\rm Li}(x)+O_{q}\left(\frac{x(\log \log x)^4}{\log^3 x}\right),
\end{equation}
where the subscript $q$ indicates that the implied constant may depend on $q$ and
$Li(x)=\int_2^x dt/\log t$ denotes the logarithmic integral.

Let ${\cal G}_q$ be the subsemigroup of the natural numbers generated by the
primes in ${\cal P}_q$. Let $HSD(x)$ count the number of non-isomorphic abelian
groups of order $n$ with $(n,q)=1$ and $n\le x$ for
which an ideal code of $F[G^*]$ whose
extension is Hermitian self-dual exists. Then by Theorem \ref{Thm concon2} we have that
$$HSD(x)=\sum_{n\le x,~n\in {\cal G}_q}a(n),$$
where $a(n)$ denotes the number of non-isomorphic abelian groups 
having $n$ elements. 

Thus we are naturally led to study the behaviour of $a(n)$
on subsemigroups $\cal G$ of the natural numbers. For our purposes it is enough 
to restrict to subsemigroups $\cal G$ that are
generated by a set $\cal P$ of primes satisfying
\begin{equation}
\label{hallo}
{\cal P}(x)=\tau {\rm Li}(x)+E_{\cal P}(x),
\end{equation}
where $0<\tau<1$ and the error term $E_{\cal P}(x)$ is small enough.

Although the literature on $a(n)$ is quite extensive, the latter problem does not seem
to have been studied before. 
Before delving into it, we recall some relevant
facts on the behaviour of $a(n)$.

\subsection{Counting non-isomorphic abelian groups}
It is easy to see that $a(n)$ is a multiplicative function with the property
that $a(p^k)=P(k)$ for every prime $p$ and every integer $k\geq 1,$ where $P(k)$
denotes the number of unrestricted partitions of $k.$ Thus $a(p^k)$ does not
depend on $p$ but only on $k,$ so that $a(n)$ is a ``prime independent''
multiplicative function.

An analytic approach to $a(n)$ is based on the fact that the Dirichlet series
associated with this function may be written as products of the Riemann zeta
function, which is defined for $\Re(s)>1$ as 
$\zeta (s)= \sum_{n=1}^{\infty} n^{-s}= \prod_p {(1-p^{-s})}^{-1}$ and otherwise
by analytic continuation.
Using the well-known identity
$$\sum_{k=0}^{\infty} P(k)x^k= \prod_{m=1}^{\infty} {1\over {{1-x}^m}},~|x|
<1,$$
one finds that, for $\Re(s)>1,$
$$\sum_{k=0}^{\infty}{{a(p^k)}\over {p^{ks}}}=\sum_{k=0}^{\infty} {{P(k)}\over
{p^{ks}}}= \prod_{m=1}^{\infty} {1\over {1-{1\over {p^{ms}}}}},$$
and thus, using the multiplicativity of $a(n),$
$$\sum_{n=1}^{\infty} {{a(n)}\over {n^s}}= \prod_p \sum_{k=0}^{\infty}
{{a(p^k)}\over {p^{ks}}}= \prod_p \prod_{m=1}^{\infty} {1\over
{1-{1\over{p^{ms}}}}}= \prod_{m=1}^{\infty} \zeta (ms).$$
Using the standard results from tauberian theory, one obtains
$$\sum_{n\leq x} a(n) \sim x \prod_{m=2}^{\infty} \zeta (m), ~ x\to \infty,$$
from this.
By much more refined methods, it can be shown that
$$\sum_{n\leq x} a(n) = \sum_{m=1}^3 c_m x^{1/m}+E(x),~c_m= \prod_{\substack{k=1
\\ k\neq m}}^{\infty}
\zeta \left({k\over m}\right),$$
where the estimates for the error term $E(x)$ have a long history
of improvements, with the best result to date being due
to Robert and Sargos \cite{RS}, who  proved that $|E(x)|\ll x^{{1/4}+\epsilon}$.
Furthermore one has, see \cite[p. 274]{F}, $c_1=2.2948565916\cdots$, 
$c_2=-14.6475663016\cdots$ and $c_3=118.6924619727\cdots$.

Thus on average $a(n)$ is constant (namely 
about $2.29$). Individual values, however, might get large.
In this direction Kr\"atzel \cite{K} proved that
\begin{equation}
\label{kra}
\lim_{n\rightarrow \infty}\sup \log (a(n)){\log \log n\over \log n}={\log 5\over 4},
\end{equation}
which implies that $a(n)\ll n^{\epsilon}$ for every $\epsilon>0$.

Ivi\'c \cite{I1} has pointed out that $C(x)$, the number of {\it distinct} values assumed
by $a(n)$ for $n\le x$, satisfies the bound
\begin{equation}
\label{cx}
C(x)\le \exp((1+o(1))2\pi\sqrt{\log x/3\log \log x}).
\end{equation}
The reason for this (see \cite[pp.  130-131]{I1}) is that there are 
$$\exp((1+o(1))2\pi\sqrt{\log x/3\log \log x})$$
integers $n\le x$ of the form
\begin{equation}
\label{form}
n=2^{a_2}3^{a_3}\cdots p^{a_p},~a_2\ge a_3\ge \cdots \ge a_p\ge 1,
\end{equation}
which is a classical result of Hardy and Ramanujan \cite[pp. 245-261]{R}.
Suppose that $a(n)$ is counted by $C(x)$, and let
$$n=p_{1,1}^{b_1}\cdots p_{1,k}^{b_k},~b_1\ge b_2\ge \cdots \ge b_k\ge 1,$$
be the canonical decomposition of $n$. Then if $m=2^{b_1}3^{b_2}\cdots p_k^{b_k}$, we
have $m\le n$ and $a(m)=P(b_1)\cdots P(b_k)=a(n)$. Therefore $C(x)$ does not
exceed the number of $n\le x$ having the form (\ref{form}) and hence
inequality (\ref{cx}) holds.

Note that if $f$ is {\it any} prime independent function, then the number of distinct
 values assumed by it for $n\le x$ satisfies the same
 upperbound as in (\ref{cx}).

\subsection{Summing $a(n)$ over $\cal G$}
Let $\chi_{\cal G}$ be the characteristic
function of $\cal G$, i.e.,
$$\chi_{\cal G}(n)=\begin{cases}
1 & {\rm if~}n{\rm ~is~in~}{\cal G};\cr 
0 & {\rm otherwise}.
\end{cases}
$$ 
We consider
$$\sum_{n\le x,~n\in {\cal G}}a(n)=\sum_{n\le x}\chi_{\cal G}(n)a(n).$$
Note that $\chi_{\cal G}(n)a(n)$ is multiplicative in $n$.
\begin{thm}
\label{blup}
If {\rm (\ref{hallo})} is 
satisfied with $E_{\cal P}(x)=O(x\log^{-1-\gamma}x)$ and $0<\gamma<1$, then
$$
\sum_{n\le x,~n\in {\cal G}}a(n)=xb_{0}\log^{\tau-1}x+O_{\cal G}(x\log^{\tau-1-\gamma/2}x).
$$
If {\rm (\ref{hallo})} is 
satisfied with $E_{\cal P}(x)=O(x\log^{-2-\gamma}x)$ and $\gamma>0$.
Then
\begin{equation}
\sum_{n\le x,~n\in {\cal G}}a(n)=
x\sum_{0\le \nu<\gamma}b_{\nu}\log^{\tau-1-\nu}x+O_{\cal G}(x\log^{\tau-1-\gamma+\epsilon}x),
\end{equation}
where $b_0,b_1,\ldots$ are constants possibly depending on $\cal G$ and
$$b_0={1\over \Gamma(\tau)}\lim_{s\downarrow 1}(s-1)^{\tau}\sum_{n\in {\cal G}}{a(n)\over n^{s}}>0.$$
\end{thm}
The proof uses the following lemma, which except for the formula for $b_0$ is
taken from \cite{MC}. The formula for $b_0$ is well-known.
\begin{lem} {\rm \cite{MC}}.
\label{een}
\label{moca}
Let $f:\mathbb N_{\ge 0}\rightarrow \mathbb R_{\ge 0}$ be a multiplicative
function satisfying
\begin{equation}
\label{groei}
0\le f(p^r)\le c_1c_2^r,~c_1\ge 1,~1\le c_2<2,
\end{equation}
and
\begin{equation}
\label{groei2}
\sum_{p\le x}f(p)=\tau {\rm Li}(x)+O({x\log^{-2-\gamma}x}),
\end{equation}
where $\tau>0$ and $\gamma>0$ are fixed, then, for $\epsilon>0$,
$$
\sum_{n\le x}f(n)=x\sum_{0\le \nu<\gamma}b_{\nu}\log^{\tau-1-\nu}x+O(x\log^{\tau-1-\gamma+\epsilon}x),
$$
where $b_0={1\over \Gamma(\tau)}\lim_{s\downarrow 1}(s-1)^{\tau}\sum_{n=1}^{\infty}{f(n)n^{-s}}$.
\end{lem}

\noindent {\it Proof of Theorem} \ref{blup}. The first assertion has been proved
by Odoni \cite{O} using a tauberian remainder theorem due to Subhankulov.
\hfill \Q{\Box}

In order to prove the second assertion we apply Lemma \ref{een} 
with $f(n)=a(n)\chi_{\cal G}(n)$. The fact that 
condition (\ref{groei}) is satisfied
follows from the classical result of Hardy and Ramanujan (see \cite[p. 240]{R}), that
$a(p^r)=P(r)=(1+o(1))(4\sqrt{3}r)^{-1}e^{\pi\sqrt{2r/3}}$ as $r$ tends to infinity.
However, the much more easily proved upperbound $P(r)\le 5^{r/4}$, see \cite{K}, is already 
sufficient in order
to show that (\ref{groei}) is satisfied.
The assumption on $E_{\cal P}(x)$ ensures that condition (\ref{groei2}) is satisfied.
On invoking Lemma \ref{moca} the proof is then completed.
   
For our problem at hand we the find the following estimate:
\begin{thm}
\label{twee}
Let $HSD(x)$ count the number of non-isomorphic abelian
groups of order $n$ with $(n,q)=1$ and $n\le x$ for
which an ideal code of $F[G^*]$ whose
extension is Hermitian self-dual exists. Then
$$HSD(x)=b_0{x\over \log^{\delta(q)}x}+O_{\epsilon,q}\left({x\over \log^{\delta(q)+1-\epsilon}x}\right),$$
where $$b_0={1\over \Gamma(1-\delta(q))}\lim_{s\downarrow 1}(s-1)^{1-\delta(q)}
\sum_{n\in {\cal G}}{a(n)\over n^{s}}.$$
\end{thm}

\subsection{The connection with free arithmetical semigroups}
A much weaker form of Theorem \ref{twee} is obtained as a straightforward consequence of
Bredikhin's Theorem, which is a basic result in the theory of {\it free arithmetical
semigroups}.

Let $G$ be a commutative semigroup with identity element 1, relative to a multplication operation
denoted by juxtaposition. Suppose that $G$ has a finite or countably infinite subset $P$ of
generators and that $G$ is {\it free}. This means that every element $n$ in $G$ has a unique
factorisation of the form $n=\omega_1^{a_1}\cdot \omega_2^{a_2}\cdots \omega_r^{a_r}$, where
the $\omega_r$ are distinct elements of $P$, the $a_i$ are possible integers, and uniqueness is
 up to order of factors. A free semigroup will be called a free arithmetical semigroup if in addition
 there exists a homomorphism of $G$ into some multiplicative semigroup ${\overline G}$ consisting
 of real numbers such that for every $x>0$, $G$ contains only finitely many elements $n$ with
 $|n|\le x$, where $|n|$ denotes the image (or norm) of the element $n$ of $G$ under the
 homomorphism $|.|$. (In the older literature the generators of ${\overline G}$ are called
 Beurling's generalized primes.) Bredikhin's theorem, for
 a proof see e.g. \cite[pp. 92-99]{Po}, then reads as follows:
\begin{thm} (Bredikhin.)
If $G$ is a free arithmetical semigroup such that
\begin{equation}
\label{weereen}
\sum_{|\omega|\le x,~\omega\in G}1=\tau{x\over \log x}+O\left({x\over \log^{1+\gamma}x}\right),
\end{equation}
where $\tau>0$ and $\gamma>0$ are fixed, then
$$
\sum_{|n|\le x\atop n\in G}1=C_Gx\log^{\tau-1}x+O(x\log^{\tau-1}x(\log \log x)^{-\gamma_1}),$$
where $\gamma_1=\min (1,\gamma)$ and
$C_G={\Gamma(\tau)}^{-1}\lim_{s\downarrow 1}(s-1)^{\tau}\sum_{n\in G}|n|^{-s}$.
\end{thm}

Now consider the free 
arithmetical semigroup $G$ of all non-isomorphic finite abelian groups with as composition the
usual direct product operation and as norm function $|A|={\rm card}(A)$. By the fundamental
theorem on finite abelian groups, $G$ is a free arithmetical semigroup having
$C(p),~C(p^2),~C(p^3),\cdots$ as generators, where $p$ runs over all the primes and $C(n)$ is
the cyclic group of order $n$. 
Since the number of cyclic groups of prime power order whose norm is not prime having
norm $\le x$ is $O(\sqrt{x}\log x)$, by the prime number theorem
in the form $\pi(x)=x/\log x + O(x/\log^2 x)$, (\ref{weereen}) is satisfied with 
$\tau=1$ and $\gamma=1$. It then follows from Bredikhin's theorem that
$$\sum_{|n|\le x\atop n\in G}1=\sum_{n\le x}a(n)=x\prod_{m=2}^{\infty}\zeta(m)+
O\left({x\over \log\log x}\right),$$
where we have used the observations that $\sum_{n\in G}|n|^{-s}=
\sum_n a(n)n^{-s}=\prod_{m=1}^{\infty}\zeta(ms)$
and $\lim_{s\downarrow 1}(s-1)\zeta(s)=1$.

Now let $G_q$ be the free arithmetical semigroup generated by all cyclic groups of the
form $C(p),C(p^2),C(p^3),\ldots$, with ord$_p(q)$ is odd or ord$_p(q^2)$ is even. Then
similarly using Bredikhin's theorem we obtain the result in Theorem \ref{twee} with the much weaker
error term $O_q(x\log^{-\delta(q)}x (\log \log x)^{-1})$.

\subsection{The maximal order of $a(n)$ on $\cal G$}
In this section we indicate what Kr\"atzel's result (\ref{kra}) looks like when one
considers the maximal order of $a(n)$ on the subsemigroup $\cal G$.
\begin{thm}
\label{drie}
Let $A=A(n)$ be the smallest integer such that $$\sum_{p\in {\cal P},~p\le A}\log p\ge 
{(\log n)}/4.$$
Then as $n$ tends to infinity and runs through the elements of $\cal G$, the estimate
$$\log a(n) \le {\cal P}(A)\log 5+O({\cal P}(A^{\theta})\log A),$$
holds with $\theta=\log(121)/\log(125)<0.994$, and there are infinitely many integers
$n$, $n\in G$, for which one has $\log a(n) ={\cal P}(A)\log 5$.
\end{thm}
{\it Proof}. Completely similar to that of the
(only) theorem in Schwarz and Wirsing \cite{SW}, who proved
this result in case $\cal P$ is the full set of primes. In their proof one merely
intersects every range of primes that occurs with $\cal P$.
\hfill \Q{\Box}

\vskip .1in
\noindent {\tt Remark}. The implicit constant in the order term can taken to be
$${2\pi^2\over 3\log 5\cdot \log 2}=5.898\cdots$$
\begin{thm}
\label{klaarhoor}
If ${\cal P}(x)\sim \tau x/\log x$ as $x$ tends to infinity, then
$$\lim_{n\rightarrow \infty}\sup_{n\in {\cal G}} 
\log (a(n)){\log \log n\over \log n}={\log 5\over 4}.$$
\end{thm}
{\it Proof}. By a standard argument in elementary
number theory it follows that if ${\cal P}(x)\sim \tau x/\log x$, 
then $\sum_{p\in {\cal P},~p\le x}\log p\sim \tau x$, $A(n)\sim (\log n)/(4\tau)$ and
${\cal P}(A)\sim \log n/(4\log \log n)$. On invoking Theorem \ref{drie} the
result then follows.
\hfill \Q{\Box}

\vskip .1in
\noindent {\tt Remark}. Let $p_1,p_2,\cdots$ denote the consecutive primes in $\cal P$.
Let $n_r=\prod_{i=1}^r p_i^4$. Suppose that ${\cal P}(x)\sim \tau x/\log x$ as $x$ tends
to infinity. We leave it as an exercise to the reader to show that
$$\lim_{r\rightarrow \infty}\log(a(n_r)){\log \log n_r\over \log n_r}={\log 5\over 4}.$$

\noindent {\tt Remark}. It is rather surprising that in Theorem \ref{klaarhoor} the
estimate does not depend on $\tau$. A similar situation arises if one compares the
maximal order of $\log d(n)$ with that of $\log r(n)$, where $d(n)$ denotes the number
of divisors of $n$ and $r(n)$ the number of way $n$ can be written as a sum of two
squares. Jacobi proved that $r(n)=4\{d_1(n)-d_3(n)\}$, where $d_1(n)$ and $d_3(n)$ denote
the number of the divisor of $n$ of the form $4k+1$ and $4k+3$, respectively. Thus
$r(n)$ counts (crudely) the divisors of $n$ made up of prime $\equiv 1({\rm mod~}4)$.
These primes have density $1/2$ amongst all primes, but nevertheless the maximal
orders of $\log d(n)$ and $\log r(n)$ are the same. Namely, we have
$$\lim_{n\rightarrow \infty}\sup \log (d(n)){\log \log n\over \log n}={\log 2},~
\lim_{n\rightarrow \infty}\sup \log (r(n)){\log \log n\over \log n}={\log 2}.$$
For further details see e.g. Nicolas \cite{Nic}. The maximal order for $\log d(n)$
was first determined by S. Wigert in 1907. Hardy and Wright \cite[Theorem 338]{HW} 
erroneously give $(\log 2)/2$ instead of $\log 2$ in the result for $\log r(n)$.

\subsection{Counting distinct values assumed by $a(n)$ on $\cal G$}
Let $C_{\cal G}(x)$ denote the number of distinct values assumed by $a(n)$
with $n\in {\cal G}$ and $n\le x$.
\begin{thm}
Let $p_0$ be the smallest prime in $\cal P$.
Suppose that there are positive constants $c_3$ and $c_4$ such that, for $x\ge p_0$,
$$c_3x< \sum_{p\in {\cal P},~p\le x}\log p<c_4x,$$
then 
$$\log C_{\cal G}(x)\sim \log C(x)\sim (1+o(1))2\pi\sqrt{\log x/3\log \log x},$$
as $x$ tends to infinity.
\end{thm}
{\it Proof}. Very similar to that given in \cite[pp. 245-261]{R}. Instead of
defining $l_n$ to be the product of the first $n$ consecutive primes, we define
it to be the product of the first $n$ consecutive primes in $\cal P$. Then instead
of (3.23) we find
$\phi(s)>c_1\int_{p_0}^{\infty}e^{-c_1sx}dx/\log x+O(1)$ and instead
of (3.24) we find $\phi(s)<c_2\int_{p_0}^{\infty}e^{-c_2sx}dx/\log x+O(1)$.
This, through Lemma 3.4, then leads to the same asymptotic for $\phi(s)$
as in the paper of Hardy and Ramanujan. This then results in the same
asymptotic for $C_{\cal G}(x)$ as that for $C(x)$.
\hfill \Q{\Box}

\vskip .3in
\noindent {\bf Acknowledgements}. The first author gratefully acknowledges financial support from the 
University of the Philippines and from the Philippine Council for Advanced Science and Technology Research 
and Development through the Department of Science and Technology.

The second author would like to thank Alexander Ivi\'c for pointing out reference \cite{RS} to him.


{\small
\begin{thebibliography}{DGS}


\bibitem{intro in handbook} R.A.~Brualdi, W.C.~Huffman, V.S.~Pless, \emph{An Introduction to Algebraic 
Codes}, in \emph{Handbook of Coding Theory}, V.S.~Pless \& W.C.~Huffman (Editors), Elsevier Science, 
Amsterdam (1998), pp. 3-139.

\bibitem{DMS} L.~Dicuangco, P.~Moree, P.~Sol\'e, \emph{The Lengths of Hermitian Self-dual Extended Duadic 
Codes}. Preprint (2005), arXiv:math.CO/0511295, submitted.

\bibitem{cunsheng} C.~Ding, D.R.~Kohel, S.~Ling, \emph{Split Group Codes}. IEEE Transactions on 
Information Theory, Vol. {\bf IT-46} (2000), pp. 485-495.

\bibitem{F} S.R.~Finch, \emph{Mathematical constants}. 
Encyclopedia of Mathematics and its Applications {\bf 94}, Cambridge University Press, 
Cambridge, 2003.

\bibitem{HW}  G.H.~Hardy, E.M.~Wright, 
\emph{An introduction to the theory of numbers}. Fifth edition. The Clarendon Press, Oxford University Press, 
New York, 1979

\bibitem{fundamentals} W.C.~Huffman, V.~Pless, \emph{Fundamentals of Error-Correcting Codes}. Cambridge 
University Press, 2003.

\bibitem{hungerford} T.W.~Hungerford, \emph{Algebra}. Springer-Verlag, New York, 1974.

\bibitem{IR} K.~Ireland and M.~Rosen, \emph{A Classical Introduction to Modern Number
Theory (Second Edition)}. Graduate Texts in Mathematics {\bf 84},
Springer-Verlag, New York, 1990.

\bibitem{I1} A.~Ivi\'c, \emph{On the number of abelian groups of a given order and on certain related 
multiplicative functions}. J. Number Theory, {\bf 16}  (1983), pp. 119-137. 

\bibitem{K} E.~Kr\"atzel, \emph{Die maximale Ordnung der Anzahl der wesentlich verschiedenen 
abelschen Gruppen $n$-ter Ordnung}. Quart. J. Math. Oxford Ser, {\bf 21} (1970), pp. 273-275. 
 
\bibitem{leon et al} J.S.~Leon, J.M.~Masley, V.~Pless, \emph{Duadic Codes}. IEEE Transactions on 
Information Theory, Vol. {\bf IT-30}, No. 5 (1984), pp. 709-714.

\bibitem{SD over GF4} F.J.~MacWilliams, A.M.~Odlyzko, N.J.A.~Sloane, H.N.~Ward, \emph{Self-Dual Codes 
over GF(4)}. Journal of Combinatorial Theory, Series {\bf A 25} (1978), pp. 288-318.

\bibitem{theory} F.J.~MacWilliams, N.J.A.~Sloane, \emph{The Theory of Error-Correcting Codes}. North 
Holland Publishing Company, Amsterdam (1983).

\bibitem{conchita} C.~Mart\'{i}nez-P\'{e}rez, W.~Willems, \emph{Self-Dual Extended Cyclic Codes}. Preprint 
(2004), http://www.math.uni-magdeburg.de/preprints/shadows/04-30report.html.

\bibitem{moree} P.~Moree, \emph{On the divisors of $a^{k} + b^{k}$}. Acta Arithmetica {\bf 80}, No. 3 
(1997), pp. 197-212.

\bibitem{Moree2} P.~Moree, \emph{On primes $p$ for which $d$ divides ord$_p(g)$}.
Funct. Approx. Comment. Math. {\bf 33} (2005), pp. 85-95.

\bibitem{MC} P.~Moree and J.~Cazaran, \emph{On a claim of Ramanujan in his first letter to Hardy}. 
Exposition. Math., {\bf 17} (1999), pp. 289-311.

\bibitem{Nic} J.-L.~Nicolas, \emph{On highly composite numbers}, in \emph{Ramanujan revisited} 
(Urbana-Champaign, Ill., 1987), Academic Press, Boston, MA, 1988, pp. 215--244.

\bibitem{O}  R.W.K.~Odoni, \emph{A problem of Rankin on sums of powers of cusp-form coefficients}. 
J. London Math. Soc. (2), {\bf  44}  (1991), pp. 203-217. 

\bibitem{qcodes} V.~Pless, \emph{Q-Codes}. Journal of Combinatorial Theory, Series {\bf A 43} (1986), pp. 
258-276.

\bibitem{onweights} V.~Pless, J.M.~Masley, J.S.~Leon, \emph{On Weights in Duadic Codes}. Journal of 
Combinatorial Theory, Series {\bf A 44} (1987).

\bibitem{Po} A.G.~Postnikov, \emph{Introduction to analytic number theory}. Translations of Mathematical 
Monographs {\bf 68}, AMS, Providence, RI, 1988. 

\bibitem{R} S.~Ramanujan, \emph{Collected papers}. Chelsea, New York, 1962.

\bibitem{rushanan} J.J.~Rushanan, \emph{Duadic Codes and Difference Sets}. Journal of Combinatorial 
Theory, Series {\bf A 57} (1991), pp. 254-261.

\bibitem{RS}  O.~Robert and P.~Sargos, \emph{Three-dimensional exponential sums with monomials}. J. Reine 
Angew. Math. {\bf 591} (2006), pp. 1-20.

\bibitem{SW} W.~Schwarz and E.~Wirsing, \emph{The maximal number of non-isomorphic abelian groups of order 
$n$}. Arch. Math. (Basel),  {\bf 24}  (1973), pp. 59-62.

\bibitem{smid} M.~Smid, \emph{Duadic Codes}. IEEE Transactions on Information Theory, Vol. {\bf IT-33} No. 
3 (1987), pp. 432-433.

\bibitem{smidthesis} M.~Smid, \emph{On Duadic Codes}. Master's Thesis, Eindhoven University of Technology, 
The Netherlands, 1986.

\bibitem{ward} H.N.~Ward, \emph{Quadratic Residue Codes and Divisibility}, in \emph{Handbook of Coding 
Theory}, V.S.~Pless \& W.C.~Huffman (Editors), Elsevier Science, Amsterdam (1998), pp. 827-870.

\bibitem{partitions} H.N.~Ward, L.~Zhu, \emph{Existence of Abelian Group Code Partitions}. Journal of 
Combinatorial Theory Series {\bf A 67} (1994), pp. 276-281.

\end{thebibliography}}
\end{document}